\documentclass[12pt]{article}
\usepackage{bm, amsmath, amssymb, amsthm}

\newcommand{\eq}{\hspace*{-2mm}&=&\hspace*{-2mm}}

\newcommand{\R}{{\mathbb R}}
\newcommand{\Sf}{{\mathbb S}}
\newcommand{\rank}{\mbox{rank$\,$}}
\newtheorem{theorem}{Theorem}
\newtheorem{proposition}[theorem]{Proposition}
\newtheorem{lemma}[theorem]{Lemma}

\newtheorem{corollary}[theorem]{Corollary}

\newtheorem{example}[theorem]{Example}

\newcommand{\spa}{\mbox{span}}
\def\<{\langle}
\def\>{\rangle}
\def\d{\partial}
\def\be{\begin{equation} }
\def\ee{\end{equation} }
\newcommand\bea{\begin{eqnarray*}}
\newcommand\eea{\end{eqnarray*}}
\begin{document}

\title{Euclidean hypersurfaces with a totally geodesic\\ foliation of
codimension one}

\author{M. Dajczer, V. Rovenski and R. Tojeiro}
\date{}
\maketitle

\begin{abstract}

We  classify the hypersurfaces of Euclidean space that carry  a
totally geodesic foliation with complete  leaves of codimension one.
In particular, we show that rotation hypersurfaces with complete 
profiles of codimension one are characterized by their warped product structure.
The local version of the problem is also considered.
\end{abstract}


A smooth foliation ${\cal F}$ on an $n$-dimensional  Riemannian manifold $M^n$
is \emph{totally geodesic} if all leaves of ${\cal F}$ are totally
geodesic submanifolds of $M^n$, that is, if any geodesic of
$M^n$ that is tangent to ${\cal F}$ at some point  is contained in the leaf
of ${\cal F}$ through that point.

Several authors have investigated whether
on a given Riemannian manifold $M^n$ there exists a totally geodesic foliation
of codimension one, that is, with $(n\!-\!1)$-dimensional
leaves, as well as the inverse problem of determining whether one can find a
Riemannian metric on a manifold $M^n$ with respect to which a given smooth
foliation of codimension one on $M^n$ becomes totally geodesic.
We refer to the work of Ghys \cite{ghys} for a complete solution of the
latter problem in the compact as well as some noncompact cases,
where several other references  can be found. See also \cite{rw} for further
discussion on the subject as well as an updated list of references.
\vspace{1ex}

In this paper, we address the following related  and rather basic extrinsic  problem: What
are all Euclidean hypersurfaces  $f\colon M^n\to\R^{n+1},\, n\geq 3$,
that carry  a  foliation of codimension one  with totally geodesic 
(complete or not) leaves?

\vspace{1ex}

First, we  discuss several families of solutions to this problem, starting
with some trivial ones.

\medskip

\noindent \underline{\emph{Flat hypersurfaces}}.  These are isometric 
immersions $f\colon U\to \R^{n+1}$ of open subsets $U\subset \R^n$, 
which admit foliations  by (open subsets of) affine hyperplanes.
Complete flat Euclidean  hypersurfaces are well known to be 
cylinders over  plane curves \cite{hn}.

\medskip

\noindent \underline{\emph{Surface-like hypersurfaces}}. For a surface 
$g\colon L^2\to\R^3$, let $\mathcal{D}_0$ be the one-dimensional distribution on
$L^2$ spanned by the tangent directions to a foliation of $L^2$ by geodesics.
 Set $M^n=L^2\times\R^{n-2}$ and define an isometric immersion
$f\colon M^n\to\R^{n+1}$ by $f=g\times id$, where $id\colon \R^{n-2}\to \R^{n-2}$
is the identity map. Then $\mathcal{D}=\mathcal{D}_0\oplus \R^{n-2}$ is clearly 
a totally geodesic distribution on $M^n$ of codimension one, whose leaves are 
complete whenever the same holds for those of $\mathcal{D}_0$.
We call $f$ a \emph{cylindrical surface-like} hypersurface.

Similar examples, but with never complete leaves,  can be constructed by 
starting with a surface in the sphere $g\colon L^2\to\mathbb{S}^3\subset \R^4$,  
with $\mathcal{D}_0$ as before, and defining  $M^n=L^2\times\R_+\times \R^{n-3}$ 
and $f=C(g)\times id$, where $C(g)\colon  L^2\times\R_+\to\R^4$, given by 
$C(g)(x, t)=tg(x)$, is the cone over $g$ in $\R^4$.
Then, the distribution $\mathcal{D}=\mathcal{D}_0\oplus \R\oplus \R^{n-3}$ on 
$M^n$ is again totally geodesic of rank $n-1$, that is, $\dim \mathcal{D}(x)=n-1$ for any $x\in M^n$. In this case, we say that $f$ 
is a \emph{conical surface-like} hypersurface.
\medskip

\noindent \underline{\emph{Ruled hypersurfaces}}. Nonflat ruled hypersurfaces
$f\colon  M^n\to \R^{n+1}$ carry a  smooth  foliation of  codimension one by
(open subsets of) affine subspaces of $\R^{n+1}$, the rulings of $f$.
Thus, the rulings are totally geodesic in $\R^{n+1}$, hence also in $M^n$.
For complete examples see Example~\ref{example}  below.
\medskip

\noindent\underline{\emph{Partial tubes}}. Let $\gamma\colon I\subset\R\to\R^{n+1}$
be a unit speed curve. Consider a hypersurface $N^{n-1}$ of the (affine) normal
space  to $\gamma$ at some  point 
and parallel transport $N^{n-1}$ along $\gamma$ with respect to the normal connection.
Then, if $N^{n-1}$ lies in a suitable open subset of that normal space
(as described in the next section), this generates an $n$-dimensional hypersurface
$M^n$ of $\R^{n+1}$, called the \emph{partial tube} over $\gamma$ with fiber $N^{n-1}$.
It turns out that the parallel translates of
$N^{n-1}$ give rise to a  totally geodesic foliation of $M^n$ of codimension one, 
whose leaves are complete if so is $N^{n-1}$.
\medskip

 Partial tubes were introduced  in \cite{cd} and  \cite{cw} and will be discussed
in more detail in the next section.  We point out that, starting with a unit speed
circle $\gamma\colon  I\to \R^2\subset \R^{n+1}$, the preceding construction
yields a rotation hypersurface of $\R^{n+1}$ having $N^{n-1}$ as profile.

The classes of hypersurfaces just described are clearly not disjoint. For instance,
the class of flat hypersurfaces is precisely the intersection of the classes of
ruled hypersurfaces and partial tubes over curves. In fact, flat hypersurfaces
free of totally geodesic points correspond to  partial tubes over  curves with 
fiber a totally geodesic hypersurface
$N^{n-1}$ of a fixed normal space to the curve.  On the other hand, a surface-like
hypersurface   is also a partial tube if and only if the integral (geodesic) curves
of $\mathcal{D}_0$ are also lines of curvature of $g$, i.e., the surface $g$ is itself a
partial tube over an orthogonal trajectory of $\mathcal{D}_0$.
Moreover, a surface-like hypersurface is flat (respectively, ruled) if and only if
the surface $g$ has
index of relative nullity one (respectively, is ruled).

In view of the discussion in the preceding paragraph, it is easy to construct
examples, even complete ones, where  different types of hypersurfaces
are smoothly attached. This is illustrated by the following simple class of examples.

\begin{example}\label{example} {\em Let $c\colon\R\to\R^{n+1}$ be a
smooth curve parametrized by
arc-length $s$ with curvatures $\kappa_1,\ldots,\kappa_n$ in a Frenet frame
$e_1=c',e_2,\ldots,e_{n+1}$. Assume that $\kappa_1>0$ along $c$ and that
$\kappa_j=0$ for $s\leq 0$ and  $\kappa_j>0$ for $s>0$,   $j\geq 2$.
Then, the complete hypersurface $F\colon\R^n\to\R^{n+1}$ parametrized by
$$
 F(s,t_1,\ldots,t_{n-1})=c(s)+\sum\nolimits_{j=1}^{n-1}t_je_{j+2}
$$
is surface-like for $s\leq 0$ and ruled but not surface-like for $s>0$.
Of course, the ruled surface factor for $s\leq 0$ can be deformed or
replaced by a non-ruled one. Notice that if $\kappa_j>0$ on all of $\R$
then $F$ is  ruled and  complete.
}\end{example}

 Under an assumption of  global nature, our main result shows that
partial tubes over curves  and ruled hypersurfaces are the only examples 
with complete leaves. Given a hypersurface $f\colon  M^n\to\R^{n+1},\,n\geq 3$, 
we say that $f(M)$ contains a \emph{surface-like strip} if
there exists an open  subset $U\subset M^n$
isometric to a product $L^2\times\R^{n-2}$ where $f$ splits as $f=g\times id$,
with $g\colon L^2\to \R^3$  an isometric immersion and
$id\colon \R^{n-2}\to\R^{n-2}$  the identity map.

\begin{theorem} \label{thm:main2}
Let $f\colon  M^n\to\R^{n+1},\,n\geq 3$, be an isometric immersion of a~nowhere
flat connected Riemannian manifold  that carries a totally geodesic
foliation of codimension one with complete leaves. If $f(M)$ does not contain any
surface-like strip then  it is either ruled or
a partial tube over a curve.
\end{theorem}

If $f\colon M^n\to\R^{n+1}$  is  an isometric immersion of a
Riemannian manifold with positive sectional curvatures, then neither $f$ can be 
ruled on any open subset nor $f(M)$ can contain  any surface-like strip. 
Thus we obtain  the following  immediate consequence of the preceding result.

\begin{corollary} \label{cor:ricci}
Let $f\colon M^n\to\R^{n+1},\, n\geq 3$, be an isometric immersion of a
Riemannian manifold with positive sectional curvatures that carries a totally 
geodesic foliation of codimension one with complete leaves.  Then  $f(M)$ 
is  a partial tube over a curve.
\end{corollary}

A Riemannian manifold $M^n$ that carries a totally geodesic foliation
of codimension one is locally (globally, if $M^n$ is simply connected and the
leaves of the foliation are complete)  isometric to a product manifold
$N^{n-1}\times \R$, with a \textit{twisted} product metric $d\sigma^2+\rho^2dt^2$,
where $d\sigma^2$ is a fixed metric on $N^{n-1}$ and 
$\rho\in C^\infty(N^{n-1}\times\R)$ (see \cite[Theorem~1]{pr}).
Thus, an equivalent statement of  Theorem~\ref{thm:main2} is that an isometric 
immersion $f\colon M^n\to\R^{n+1},\,n\geq 3$, of a~twisted product manifold
$M^n=N^{n-1}\times_{\rho} I$, where $I\subset \R$ is an open interval
and $N^{n-1}$ is a complete manifold free of flat points, is either ruled
or a partial tube over a curve, as soon as $f(M)$ does not contain any
surface-like strip.
\vspace{1ex}

For isometric immersions $f\colon M^n\to\R^{n+1},\, n\geq 3$, of a \emph{warped} 
product manifold $M^n=N^{n-1}\times_{\rho} I$, in which case $\rho$  depends only on 
$N^{n-1}$,  we prove the following result.

\begin{theorem} \label{thm:main3}
Let $f\colon  M^n\to \R^{n+1},\, n\geq 3$, be an isometric immersion of a warped product
connected Riemannian manifold $M^n=N^{n-1}\times_{\rho} I$ where
$N^{n-1}$ is a complete manifold free of flat points,
$\rho\in C^\infty(N)$ and $I\subset \R$ is an open interval.
If $f(M)$ does not contain any
surface-like strip then it is a rotation hypersurface having $N^{n-1}$ as profile.
\end{theorem}

Note that the assumption that $f(M)$ does not contain any
surface-like strip here means that $N^{n-1}$ does not contain any open subset $U$
isometric to $\R^{n-2}\times J$, where $J\subset \R$ is an open interval,
such that $f|_{U\times I}$ splits as $f=id\times f_1$, with
$f_1\colon J\times_{\rho} I\to \R^3$  an isometric immersion.
\vspace{1ex}

The case in which $M^n$ is assumed to be compact in Theorem \ref{thm:main3} was already
considered in \cite{mt}. Isometric immersions $f\colon  M^n\to \R^{n+1},\,n\geq 3$,
of a warped product connected Riemannian manifold free of flat points
$M^n=N^{n-k}\times_{\rho} L^k$, with $k\geq 2$, were shown in \cite{dt2}  to be, even
locally,  either rotation hypersurfaces, products with $\R^k$ of
hypersurfaces of $\R^{n-k+1}$  or products with $\R^{k-1}$ of
cones over  hypersurfaces of $\Sf^{n-k+1}\subset \R^{n-k+2}$.
\vspace{1ex}

Next, we consider the  local version of the problem stated in the beginning of 
the introduction. We prove that exactly one further class of examples may occur.

\begin{theorem} \label{thm:mainlocal}
Let $f\colon  M^n\to\R^{n+1},\,n\geq 3$, be an isometric immersion of~a Riemannian 
manifold that carries a totally geodesic  foliation of codimension one. Then, 
there exists an open dense subset of $M^n$ where  $f$ is locally either surface-like, 
ruled, a partial tube over a curve or an envelope of a~one-parameter family of flat 
hypersurfaces.
\end{theorem}

A hypersurface $f\colon M^n\to\R^{n+1}$  is called the \emph{envelope}
of a one-parameter family of hypersurfaces $F_t\colon  M_t^n\to\R^{n+1}$ if
there exists an integrable smooth distribution $\mathcal{D}$ of rank $n-1$ on $M^n$
for each leaf $\sigma_t$ of which one has an embedding 
$j_t\colon\sigma_t\to M^n_t$ such that  $F_t\circ j_t=f|_{\sigma_t}$ and 
${F_t}_*T_{j_t(x)}M_t=f_*T_xM$ for any $x\in \sigma_t$.

Geometrically, the hypersurface $F_t$ is tangent to $f$ along the  leaf $\sigma_t$ of 
$\mathcal{D}$, called the \emph{characteristic} of the one-parameter family$\{F_t\}$ at
level $t$. If the one-parameter family of hypersurfaces $F_t$ is locally defined by
the equation $G(t,x)=0$  where $x=(x_1, \ldots, x_{n+1})$,
then the envelope of the family  is locally given by
$$
\left\{\begin{array}{l}G(t,x)=0\vspace{.5ex}\\
G_t(t,x)=0,\end{array}\right.
$$
where the subscript denotes partial derivative with respect to $t$.
The characteristic at level $t=t_0$ is then the set of solutions of the
preceding system for $t=t_0$.
In Theorem \ref{thm:mainlocal}, the leaves of the  totally geodesic
foliation of $M^n$ are precisely the characteristics of the   one-parameter
family of flat hypersurfaces that envelope $f$.

\section{Partial tubes}

We first recall the precise definition of a partial tube, and then state a
result from \cite{to} that implies that partial tubes over curves are precisely
the solutions to our problem for which the orthogonal trajectories to the 
totally geodesic foliation  are lines of curvature of the hypersurface.
\vspace{1,5ex}

Let $\gamma\colon I\subset\R\to\R^N$ be a unit speed curve and let 
$\{\xi_1,\ldots,\xi_m\}$ be an orthonormal set of parallel normal vector fields along 
$\gamma$. Hence, the vector subbundle $E=\spa \{\xi_1, \ldots,\xi_m\}$ of the normal bundle
of $\gamma$ is parallel and flat. Then the map $\phi\colon I\times\R^m\to E$
given by
$$
\phi_{s}(y)=\phi(s,y)=\sum\nolimits_{i=1}^my_i\xi_i(s)
$$
for $s\in I$ and $y=(y_1, \ldots, y_m)\in \R^m$, is a parallel vector bundle isometry. 
Let $f_0\colon N^{n-1}\to\R^m$ be a substantial isometric immersion, i.e., 
an immersion whose codimension cannot be reduced.
Denote $M^n=N^{n-1}\times I$ and define a map $f\colon M^n\to\R^N$ by
$$
f(p,s)=\gamma(s)+\phi_{s}(f_0(p)).
$$
One can check that $f$ is an immersion whenever  
$f_0(N)\subset \Omega(\gamma;\phi)$, where
$$
\Omega(\gamma;\phi)=\{Y\in\R^m:\<\gamma''(s),\phi_{s}(Y)\>\neq 1 
\;\mbox{for any}\; s\in I\}.
$$
In this case, we say that $f(M)$  is  the \emph{partial tube over $\gamma$ with fiber $f_0$}.
Endowing $M^n$ with the induced metric, the distribution on $M^n$
given by the tangent spaces to the first factor is totally geodesic. Moreover, the second
fundamental form of $f$ satisfies $\alpha_f(X,\d/\d s)=0$ for any $X\in TN$, where $\d/\d s$ is a 
unit vector field tangent to  the factor $I$.

If $M^n$ is not simply connected and $\pi\colon \tilde M^n\to M^n$ is its universal covering, 
then the map $\tilde f=f\circ \pi$ satisfies $\tilde f(\tilde M)=f(M)$. Therefore, in proving
that $f(M)$ is a partial tube there is no loss of generality in assuming that $M^n$ is simply 
connected. \medskip

The next result is a direct consequence of Theorem $16$ in \cite{to}. 

\begin{theorem}\label{pt2} Let $f\colon  M^n\to\R^N$ be an isometric immersion of a twisted 
product $M^n =N^{n-1}\times_{\rho} I$, where $I\subset \R$ is an open interval and 
$\rho\in C^{\infty}(M)$. Assume that the second fundamental form of $f$  satisfies
\be\label{eq:sff}
\alpha_f( X,\d/\d s)=0\;\;\mbox{for any}\;\; X\in TN.
\ee
 Then $f(M)$ is  a partial tube over a curve.
\end{theorem}

In view of the discussion after Corollary \ref{cor:ricci}, this yields the following
result.

\begin{corollary}\label{pt} Let $f\colon  M^n\to\R^{n+1}$ be an isometric immersion carrying a
smooth totally geodesic foliation  of rank $n-1$  whose orthogonal trajectories  are lines of 
curvature of $f$. 
Then $f(M)$ is locally a partial tube over a curve.
\end{corollary}

\section{Hypersurfaces with a curvature invariant distribution}

In this section we prove a  preliminary result on oriented Euclidean hypersurfaces
 that  carry a curvature invariant distribution.
\vspace{1,5ex}

That a  smooth  distribution $\mathcal{D}$ on a Riemannian manifold
$M^n$ is  \emph{curvature invariant} means that
\be\label{eq-01}
R(X,Y)Z \in \mathcal{D}\;\;\mbox{for all}\;\; X,Y,Z\in \mathcal{D},
\ee
where $R$ denotes the curvature tensor of $M^n$.  By $\mathcal{D}$ being 
\emph{totally geodesic} we mean that 
$$
\nabla_X Y\in \mathcal{D}\;\;\mbox{for all}\;\; X,Y\in \mathcal{D}.
$$
In particular, any totally geodesic distribution $\mathcal{D}$ in the above sense 
is curvature invariant and integrable, and its leaves are totally geodesic submanifolds 
of $M^n$.  For results on curvature invariant distributions we refer to \cite{T} 
and the references therein.\vspace{1ex}

For an oriented Euclidean hypersurface, we always denote by $A$ its shape operator 
with respect to a globally defined smooth unit normal vector field, and by 
$\Delta=\ker A$ its \emph{relative nullity} distribution.

\begin{proposition}\label{prop:curv}
Let $f\colon M^{n}\to\R^{n+1}$ be an oriented hypersurface carrying a
curvature invariant distribution  $\mathcal{D}$ of rank $k>1$.
Then one of the following possibilities holds pointwise:
\begin{itemize}
\item[(i)] $A(\mathcal{D})\subset \mathcal{D}^\perp$,
\item[(ii)] $A(\mathcal{D})\subset \mathcal{D}$,
\item[(iii)] $\rank \mathcal{D}\cap \Delta=\rank \mathcal{D}-1$.
\end{itemize}
\end{proposition}

\proof By the Gauss equation
$$
R(X,Y)Z=\<AY,Z\>AX-\<AX,Z\>AY
$$
we have that (\ref{eq-01}) on a hypersurface is equivalent to
\be\label{eq-02}
\<AY,Z\>AX-\<AX,Z\>AY \in \mathcal{D}\;\;\mbox{for all}\;\; X,Y,Z\in \mathcal{D}.
\ee
At a point, let $X_1,\ldots X_k$ be a local orthonormal base of $\mathcal{D}$ such that
$$
A X_j=\lambda_j X_j +V_j
$$
where $V_j\in \mathcal{D}^\perp$ and $\lambda_j\in\R$. If  $\lambda_j=0$ for $1\le j\le k$
then $(i)$ holds. Otherwise,  we may assume that $\lambda_1\ne 0$.
Applying (\ref{eq-02}) for $Y=Z=X_1$ and $X=X_j$, $j\ge 2$, we obtain
$$
\lambda_1(\lambda_jX_j+V_j)\in \mathcal{D},
$$
hence $V_j=0$ for $2\leq j\leq k$. On the other hand,  we have from
(\ref{eq-02}) for $Y=Z=X_j$, $ j\ge 2$, and $X=X_1$ that
$$
\lambda_j(\lambda_1X_1+V_1)\in \mathcal{D}.
$$
We conclude that  either $V_1=0$ or
$\lambda_j=0$ for $2\leq j\leq k$, which correspond to cases $(ii)$ and $(iii)$,
respectively.\qed

\section{The global case}

We first construct a suitable tangent orthonormal frame on hypersurfaces that 
carry a totally geodesic distribution of codimension one and satisfy
condition $(iii)$ in Proposition \ref{prop:curv}.

\begin{lemma}\label{le:frame} Let $f\colon  M^n\to\R^{n+1}$ be an isometric 
immersion of a Riemannian manifold  that carries a totally geodesic distribution 
$\mathcal{D}$ of rank $n-1$ such  that condition $(iii)$ in 
Proposition \ref{prop:curv} holds everywhere.
Then, there exist locally a smooth orthonormal tangent frame
$\{Y,X,T^1,\ldots,T^{n-2}\}$ and smooth functions $\beta,\mu,\rho$ and $\lambda_j$,
$1\leq j\leq n-2$, such that
\be\label{shape}
AY=\beta Y+\mu X,\quad AX=\mu Y+\rho X,\;\;AT^j=0,
\ee
\be\label{E-parallel}
\nabla_{T^j}{T^i} = \nabla_{T^j}X = \nabla_{T^j}Y =\nabla_XY=0
\ee
and
\be\label{eq:ntx}
\nabla_XT^j=\lambda_j X,\;\;\;1\leq i,j\leq n-2.
\ee
Moreover, condition $(i)$ (respectively, $(ii)$) in Proposition \ref{prop:curv}  
holds at $x\in M^n$ if and only if $\rho$ (respectively, $\mu$) vanishes at $x$.
\end{lemma}

\proof
Choose unit vector fields  $Y\in \mathcal{D}^\perp$ and  $X\in \mathcal{D}$ 
orthogonal to $\Delta$. Then there exist smooth functions $\beta,\mu,\rho$ such 
that the first two equations in (\ref{shape}) are satisfied. Since $\mathcal{D}$ 
is totally geodesic, the last equation in (\ref{E-parallel}) holds.

Let $\gamma\colon J\to M^n$ be the integral curve of $Y$ through an arbitrary
given point in $M^n$,  and define $\psi\colon I\times J\to M^n$ by requiring that,
for any fixed $y\in J$, the map $x\mapsto \psi_{\gamma(y)}(x)=\psi(x, \gamma(y))$
is the integral curve of $X$ such that  $\psi_{\gamma(y)}(0)=\gamma(y)$.

The normal space in $\R^{n+1}$ of the restriction $f|_\sigma$  of $f$ to a leaf 
$\sigma$ of $\mathcal{D}$ at each point $x\in \sigma$ is spanned by  $f_*Y(x)$ 
and $\eta(x)$, where $\eta$ is a smooth unit normal vector field to $f$. 
Since $\sigma$ is totally geodesic in $M^n$, the shape operator of $f|_\sigma$ 
at $x$ with respect to $f_*Y(x)$ is identically zero. 
On the other hand, since condition $(iii)$ in Proposition \ref{prop:curv} holds 
at $x$ by assumption, the shape operator of $f|_\sigma$ at $x$  with respect to 
$\eta$ has rank one. It follows from the Gauss equations for $f|_\sigma$ that 
$\sigma$ is flat.

To construct the desired smooth orthonormal tangent frame, start with any smooth 
orthonormal frame $\{T^1, \ldots, T^{n-2}\}$ spanning $\Delta$ along~$\gamma$. Then, 
for each integral curve  $x\mapsto \psi_{\gamma(y)}(x)$ of $X$, extend $\{T^1,\ldots,T^{n-2}\}$
to an orthonormal frame spanning $\Delta$ along $\psi_{\gamma(y)}$ by parallel translation
with respect to the  normal connection of $\psi_{\gamma(y)}$ as a curve in the
leaf of $\mathcal{D}$ containing  $\psi_{\gamma(y)}(I)$. Finally, parallel translate
each $T^j(\psi_{\gamma(y)}(x))$, $1\leq j\leq n-2$, along the leaf of $\Delta$
through $\psi_{\gamma(y)}(x)$.  Since $\Delta$ and $\mathcal{D}$ are  both totally
geodesic, it follows that the orthonormal frame  $\{T^1, \ldots, T^{n-2}\}$
constructed in this way satisfies the last equation in (\ref{shape}), the first
three equations in (\ref{E-parallel}) as well as~(\ref{eq:ntx}). 
The last assertion is clear.\qed

\begin{lemma}\label{surflike} Let $f\colon  M^n\to\R^{n+1}$ be an isometric immersion
of a nowhere flat  Riemannian manifold  that carries a totally geodesic distribution
$\mathcal{D}$ of rank $n-1$ with complete leaves. Let $U$ be the open subset of $M^n$ 
where neither of conditions $(i)$ or $(ii)$ in Proposition \ref{prop:curv} occur.
Then each connected component of $U$ is  isometric to a product $W=L^2\times \R^{n-2}$ 
and $f|_W=g\times id$, where $g\colon L^2\to \R^3$ is an isometric immersion.
\end{lemma}

\proof Since any
totally geodesic distribution is curvature invariant,  it follows from
Proposition \ref{prop:curv}  and the assumptions that condition $(iii)$ holds
everywhere on $U$. Let $\{Y,X,T^1,\ldots,T^{n-2}\}$ be the frame on $U$ given by
Lemma~\ref{le:frame}. Straightforward computations using the Codazzi  equations yield
\begin{subequations}
\begin{eqnarray}
\label{En-sys1-CE-c}
T^i(\rho)\eq\rho\,\<\nabla_X X, T^i\>,\\
\label{En-sys1-CE-d}
T^i(\mu) \eq \mu\,\<\nabla_X X, T^i\>,\\
\label{En-sys1-CE-e}
T^i(\mu) \eq \rho\,\<\nabla_Y X, T^i\> +\mu\,\<\nabla_Y Y, T^i\>
\end{eqnarray}
\end{subequations}
whereas the  Gauss equations give
\begin{subequations}
\begin{eqnarray}\label{En-sys1-GE}
 X(\<\nabla_Y Y,T^i\>)\eq \<\nabla_Y Y, X\>\,(\<\nabla_Y Y,T^i\>-\<\nabla_X X, T^i\>),\\
\label{En-sys1-GE-c}
 T^i(\<\nabla_X X, T^j\>)\eq\<\nabla_X X, T^i\>\,\<\nabla_X X, T^j\>,\\
\label{En-sys1-GE-d}
 T^i(\<\nabla_Y Y, T^j\>)\eq\<\nabla_Y Y, T^i\>\,\<\nabla_Y Y, T^j\>.
\end{eqnarray}
\end{subequations}

Since neither of conditions $(i)$ or $(ii)$ in Proposition \ref{prop:curv} is
satisfied at any point of $U$, we have that $\mu\rho\neq 0$ everywhere on $U$.
On the other hand, from equations (\ref{En-sys1-CE-c}) and (\ref{En-sys1-CE-d}) 
we obtain that
\be\label{vrho} 
\mu=\varphi\rho,\;\mbox{where}\;\;T^j(\varphi)=0\;\mbox{for}\;1\leq j\leq n-2.
\ee

Consider a unit speed geodesic $\gamma$ starting at a  point of $U$ and tangent to some $T^j$. 
Since the leaves of $\mathcal{D}$ are assumed to be complete,  $\gamma$ is defined at 
any value of the parameter. We claim that it remains indefinitely in $U$. 
Otherwise $\gamma$ would reach a point $y$
of the boundary of $U$. Since $y$ is the limit of a sequence
of points where either of conditions $(i)$ or $(ii)$ in Proposition \ref{prop:curv}
holds, then either $\mu$ or $\rho$ must vanish at $y$. But then both 
$\mu$ and $\rho$ vanish at $y$, in view of (\ref{vrho}). Hence $y$ is a flat 
point of $M^n$, contradicting our assumption and proving our claim.

It follows that the leaves of $\Delta$ through points of $U$ are
complete and that condition $(iii)$ remains valid along them. 
In view of (\ref{En-sys1-GE-c}) and (\ref{En-sys1-GE-d}) for $i=j$, 
 we conclude that the
functions  
\be\label{lambdaj} 
\lambda_j=-\<\nabla_XX,T^j\>\;\;\mbox{and}\;\;\theta_j
=\<\nabla_YY,T^j\>,\;\;1\leq j\leq n-2,
\ee
must be everywhere vanishing along such leaves.
In general, we obtain from  (\ref{En-sys1-CE-d}) and (\ref{En-sys1-CE-e}) that
\be\label{codazzis}
\mu(\lambda_j+\theta_j)+\rho\,\<\nabla_Y X,T^j\>
=0\;\;\mbox{for all}\;\;1\leq j\leq n-2.
\ee
Since now $\lambda_j=0=\theta_j$  and $\rho\neq 0$,
we  must have that  $\<\nabla_YX,T^j\>=0$.  In particular, this implies
that $\Delta^\perp$ is totally geodesic. Moreover, from \mbox{$T^j\in\Delta$} 
we obtain that $\tilde\nabla_Xf_*T^j=0$ and $\tilde\nabla_Y f_*T^j\in f_*\Delta$,
where $\tilde\nabla$  stands for the connection in the Euclidean ambient space.
It follows that  $f_*\Delta$  is a parallel subbundle of $f^*T\R^{n+1}$,
and the result follows by choosing $L^2$ in each connected component of $U$ 
as a maximal integral leaf of $\Delta^\perp$.\qed
\vspace{1,5ex}

\noindent {\em Proof of Theorem \ref{thm:main2}:} In view of the assumption that 
$f(M)$ does not contain any surface-like strip, it follows from Lemma \ref{surflike} 
that either of conditions $(i)$ or $(ii)$ must hold
at any point of $M^n$. Let $S_1$ (respectively, $S_2$) be the subset of $M^n$ where
condition $(i)$ (respectively, condition $(ii)$) is satisfied.  Since both $S_1$ and
$S_2$ are closed and $M^n=S_1\cup S_2$,  any point on $\partial S_1$  belongs
to $S_1\cap S_2$, and hence is a flat point. It follows from our assumption that
either $M^n=S_1$ or $M^n=S_2$. In the first case $f$ is a ruled hypersurface.
In the latter, as pointed out before the statement of Theorem \ref{pt2}, there 
is no loss of generality in assuming that $M^n$ is simply connected. In this case, 
by the assumption that the leaves of the totally geodesic foliation are complete, 
it follows from \cite[Theorem~1]{pr} that  $M^n$ is isometric to a
twisted product $N^{n-1}\times_{\rho} \R$, with $\rho\in C^{\infty}(N^{n-1}\times \R)$,
and the fact that $M^n=S_2$ means that the second fundamental form of $f$
satisfies (\ref{eq:sff}). Thus $f(M)$ is a partial tube over a curve by Theorem \ref{pt2}.
\vspace{1,5ex}\qed

\noindent {\em Proof of Theorem \ref{thm:main3}:} From the equivalent
form of Theorem \ref{thm:main2} discussed right after the statement of Corollary \ref{cor:ricci}, it follows that  $f(M)$ must be
either ruled or a partial tube over a curve.

We now show that, if  $M^n$ is a warped product, then the first possibility can not occur under
 our global assumptions. In fact, the warped product structure on $M^n$
implies that the distribution on $M^n$ given by the tangent spaces to
the fibers corresponding to the second  factor $I$ is spherical. This means that the 
integral curves of a unit vector field $Y$ tangent to $I$ are \textit{extrinsic circles} 
in $M^n$, that is, $\nabla_Y\nabla_YY$ is everywhere a multiple of $Y$, or equivalently,
$$
\<\nabla_Y\nabla_YY, X\>=0=\<\nabla_Y\nabla_YY, T^j\>,\;\;1\leq j\leq n-2,
$$
where $\{Y,X,T^1,\ldots,T^{n-2}\}$ is the frame given by Lemma \ref{le:frame}.

As in the proof of Lemma \ref{surflike}, completeness of $N^{n-1}$ implies that
the functions 
$\theta_j$  in (\ref{lambdaj}) must vanish everywhere for $1\leq j\leq n-2$. From
$$
\<\nabla_Y\nabla_YY, T^j\>=0, \;\;1\leq j\leq n-2,
$$
we obtain 
\be\label{eq:vanish}
0=-\<\nabla_YY, \nabla_YT^j\>=\<\nabla_Y Y,X \>\<\nabla_YX, T^j\>,\;\;1\leq j\leq n-2.
\ee
On the other hand, using that $f$ is ruled we obtain from the Gauss equation for $f$ that
$$
X\<\nabla_Y Y, X\>=-\mu^2 +\<\nabla_Y Y, X\>^2,
$$
hence $\<\nabla_YY, X\>$ can not vanish on any open subset of $M^n$, because
$\mu$ is nowhere vanishing. We conclude from (\ref{eq:vanish}) that
$$
\<\nabla_Y X, T^j\>=0,\;\;1\leq j\leq n-2.
$$
But, as in the proof of Lemma \ref{surflike} this implies that
$f$ is a cylindrical surface-like hypersurface, contradicting our assumption.

We conclude that $f$ can not be ruled,  hence it is a partial tube over a~curve.
Since $M^n$ is a warped product, it follows from the main result in \cite{nol}
(see also \cite[Theorem~30]{to}),
that $f$ is either a cylinder over a plane curve, the product with an Euclidean
factor $\R^{n-2}$ of a cone over a curve in $\Sf^2\subset \R^3$
or a rotation hypersurface. The first two possibilities are ruled out by our
assumptions. \qed

\section{The local case}

In this section we prove  Theorem \ref{thm:mainlocal} in the introduction.\vspace{1,5ex}

\noindent \emph{Proof of Theorem \ref{thm:mainlocal}:} 
Let  $f\colon M^n\to\R^{n+1}$,\, $n\geq 3$, be an isometric immersion of a  Riemannian
manifold  that carries a totally geodesic distribution $\mathcal{D}$ of rank $n-1$.
Denote as before by $S_1$ (respectively, $S_2$)  the subset of $M^n$ where
condition $(i)$ (respectively, condition $(ii)$) is satisfied. Then $f$ is
ruled on the interior of $S_1$ and,  by Corollary \ref{pt}, it is
locally a partial tube over a curve on the interior of $S_2$. It follows
from Proposition~\ref{prop:curv}
that condition $(iii)$ holds everywhere on  the
open subset ${\cal U}=M^n\setminus (S_1\cup S_2)$.

We now show that $f$ is locally the envelope of a one-parameter family of flat 
hypersurfaces on ${\cal U}$. More precisely, we prove that for each leaf $\sigma_t$ 
of $\mathcal{D}$ there exist a  flat hypersurface $F_t\colon  V_t^n\to\R^{n+1}$ and 
an embedding $j_t\colon \sigma_t\to V^n_t$ such that  $F_t\circ j_t=f|_{\sigma_t}$ and 
${F_t}_*T_{j_t(x)}V_t=f_*T_xM$ for any $x\in \sigma_t$.

Consider  a smooth orthonormal frame
$\{Y,X,T^1,\ldots,T^{n-2}\}$ and smooth functions $\beta,\mu,\rho$ and $\lambda_j$,
$1\leq j\leq n-2$, given locally in ${\cal U}$ by Lemma \ref{le:frame}.
Let $\sigma_t$ be a leaf of $\mathcal{D}$ on ${\cal U}$ and set $f_t=f|_{\sigma_t}$.
Then, the normal space of  $f_{t}$  at each $x\in \sigma_t$ is spanned by $f_*Y(x)$ 
and $\eta(x)$, where $\eta(x)$ is a  unit normal vector  to $f$ at $x$.
Let $\pi_t\colon \Lambda_t\to \sigma_t$ be the line subbundle of  the normal bundle of 
$f_t$ that is spanned by the (restriction to $\sigma_t$ of the) vector field
$Z=\rho f_*Y-\mu f_*X$, and define $F_t\colon V_t\to \R^{n+1}$
as the restriction of the map
$$
\lambda\in\Lambda_t\mapsto f(\pi_t(\lambda))+\lambda
$$
to a tubular neighborhood $V_t$ of the $0$-section $j_t\colon \sigma_t\to \Lambda_t$
of $\Lambda_t$ where it is an immersion. For each $\lambda\in V_t$, we have
$T_\lambda V_t=T_{\pi_t(\lambda)}\sigma_t\oplus \Lambda_t(\pi_t(\lambda))$ where 
$T_{\pi_t(\lambda)}\sigma_t$ is identified with its horizontal lift  and 
$\Lambda_t(\pi_t(\lambda))$ with the vertical subspace at $\lambda$.
We prove below that  $F_t$ defines  a flat 
hypersurface by showing  that the subspaces  ${F_t}_*T_\lambda V_t$ are constant 
along the distribution $\tilde \Delta$ of rank $n-1$ on $V_t$ given as follows. 
For each $\lambda\in V_t$, let $\Delta(\pi_t(\lambda))$ be the relative nullity subspace 
of $f_{t}$ at $\pi_t(\lambda)$ and define  $\tilde\Delta(\lambda)
=\Delta(\pi_t(\lambda))\oplus \Lambda_t(\pi_t(\lambda))$, 
where $\Delta(\pi_t(\lambda))$ is identified with its horizontal lift.   

Given $\lambda=sZ/|Z|\in \Lambda_t$ and $W\in T_{\pi_t(\lambda)}\sigma_t$,  we have
\be\label{eq:Ft*}{F_t}_*W=f_*W+\tilde \nabla_W sZ/|Z|=f_*W+sW(|Z|^{-1})Z+s|Z|^{-1}\tilde \nabla_W Z.\ee
Since
$$\<\tilde \nabla_{X} Z, \eta\>=-\<Z, \tilde \nabla_{X} \eta\>=0,$$
for $\tilde \nabla_{X} \eta=-f_*A_\eta X=-\mu f_*Y-\rho f_*X$, we obtain that 
\be\label{eq:nx}
\<{F_t}_*X, \eta\>=0.
\ee

On the other hand,  from (\ref{E-parallel}) we have that  $\nabla_{T^j}Y=0=\nabla_{T^j}X$, whereas
equations (\ref{En-sys1-CE-c}) and 
(\ref{En-sys1-CE-d}) yield
$$
\frac{T^j(\rho)}{\rho}=\frac{T^j(\mu)}{\mu}.
$$
Using that $\Delta=\spa\{T^1, \ldots, T^{n-2}\}$ is the relative nullity 
distribution of $f$ we obtain
\begin{eqnarray*}\tilde \nabla_{T^j} Z \!\!\!&=&\!\!\! T^j(\rho)f_*Y
+\rho \tilde \nabla_{T^j}f_*Y-T^j(\mu)f_*X-\mu \tilde \nabla_{T^j}f_*X\vspace{1ex}\\
\!\!\!&=&\!\!\! \frac{T^j(\mu)}{\mu}Z.
\end{eqnarray*}
In particular, 
\begin{eqnarray*}
T^j(|Z|^2)&=&2\<\tilde \nabla_{T^j} Z, Z\>\vspace{1ex}\\
&=&2\frac{T^j(\mu)|Z|^2}{\mu},
\end{eqnarray*} and hence
$$T^j(|Z|^{-1})Z+|Z|^{-1}\tilde \nabla_{T^j} Z=0.$$
It  follows from  (\ref{eq:Ft*}) that
${F_t}_*T^j=f_*T^j$ for $1\leq j \leq n-2$, 
 and thus
 \be\label{eq:ntj}{F_t}_*\Delta(\pi_t(\lambda))=f_*\Delta(\pi_t(\lambda)).\ee  
We obtain from  (\ref{eq:nx}) and (\ref{eq:ntj}), together with  the fact that ${F_t}_*\Lambda(\pi_t(\lambda))$ is spanned by $f_*Z$, 
 that $\eta(\pi_t(\lambda))$ is normal to ${F_t}_*T_\lambda V_t$
along $\tilde \Delta$. In particular,  ${F_t}_*T_{j_t(x)}V_t=f_*T_xM$ for any $x\in \sigma_t$.\vspace{1ex}\qed

We point out that the  construction of the flat hypersurface 
$F_t\colon V_t\to \R^{n+1}$ extending $f_t$ in the proof of Theorem \ref{thm:mainlocal} is a special case of the ruled extension 
of submanifolds given in \cite[Proposition~8]{dt}.

{\renewcommand{\baselinestretch}{1} \hspace*{-20ex}\begin{tabbing}
\indent \= Marcos Dajczer\\
\> IMPA \\
\> Estrada Dona Castorina, 110\\
\> 22460-320 -- Rio de Janeiro -- Brazil\\
\> marcos@impa.br\\
\end{tabbing}}

\vspace*{-5ex}

{\renewcommand{\baselinestretch}{1} \hspace*{-20ex}\begin{tabbing}
\indent \= Vladimir Rovenski \\
\> University of Haifa \\
\> Mount Carmel, \\
\> 31905 -- Haifa -- Israel\\
\> rovenski@math.haifa.ac.il
\end{tabbing}}

\vspace*{-2ex}

{\renewcommand{\baselinestretch}{1} \hspace*{-20ex}\begin{tabbing}
\indent \= Ruy Tojeiro\\
\> Departamento de Matematica, \\
\> Universidade Federal de S\~{a}o Carlos,\\
\> 13565-905 -- S\~{a}o Carlos -- Brazil\\
\> tojeiro@dm.ufscar.br
\end{tabbing}}


\begin{thebibliography}{lll}

\bibitem{cd} S. Carter and U. Dursun,  \emph{Partial tubes and Chen submanifolds},
J. Geom.  \textbf{63} (1998), 30--38.

\bibitem{cw} S. Carter and A. West, \emph{Partial tubes about immersed manifolds},
Geom. Dedicata~{\bf 54} (1995), 145--169.

\bibitem{dt2} M. Dajczer and R. Tojeiro, \emph{Isometric immersions  in codimension
two of warped products into space forms products},
Illinois J. Math. {\bf 48} (2004), 711--746.

\bibitem{dt} M. Dajczer and R. Tojeiro, \emph{Submanifolds with nonparallel first
normal bundle revisited},  Pub. Mat.  {\bf 58} (2014), 179-191. 

\bibitem{ghys} E. Ghys, \emph{Classification des feuilletages totalement
g\'eod\'esiques de codimension un}, Comm. Math. Helv. {\bf 58} (1983), 543--572.

\bibitem{hn} P. Hartman and L. Nirenberg, \emph{On spherical image maps whose
Jacobians do not change signs}, Amer. J. of Math. {\bf 81} (1959), 901--920.

\bibitem{mt} I. Moutinho and R. Tojeiro, \emph{Polar actions on compact
Euclidean hypersurfaces},
Annals Global An. Geom. {\bf 33}  (2008), 323--336.

\bibitem{nol} S. N\"olker, \emph{Isometric immersions of warped products},
Diff. Geom. Appl.~{\bf 6} (1996), 31--50.

\bibitem{pr} R. Ponge and H. Reckziegel, \emph{Twisted products in pseudo-Riemannian
geometry}, Geom. Dedicata {\bf 48} (1993), 15--25.

\bibitem{rw}  V. Rovenski and P. Walczak, ``Topics in extrinsic geometry of
codimension one foliations",  Springer briefs in Mathematics,  Springer, 2011.

\bibitem{to}   R. Tojeiro, \emph{A decomposition theorem for immersions of product
manifolds}.  To appear in Proc. Edinburgh Math. Soc. See arXiv:1306.3278.

\bibitem{T} K. Tsukada, \emph{Totally geodesic submanifolds of Riemannian manifolds and
curvature-invariant subspaces}, Kodai Math. J. {\bf 19} (1996),  395--437.
\end{thebibliography}
\end{document}